\documentclass[11pt]{amsart}
\usepackage{a4wide}
\usepackage{amsmath}
\usepackage[utf8]{inputenc}
\usepackage{amssymb}
\usepackage{amsopn}
\usepackage{epsfig}
\usepackage{amsfonts}
\usepackage{latexsym}
\usepackage{graphicx}
\usepackage{enumerate}

\newtheorem{theorem}{Theorem}[section]
\newtheorem{lemma}[theorem]{Lemma}

\newtheorem{corollary}[theorem]{Corollary}

\theoremstyle{definition}

\theoremstyle{remark}

\numberwithin{equation}{section}

\newcommand{\R}{\ensuremath{\mathbb{R}}}
\newcommand{\N}{\ensuremath{\mathbb{N}}}

\begin{document}

\title{On the sum of squares of middle-third Cantor set}

\author[Z. Wang]{Zhiqiang Wang}
\address[Z. Wang]{School of Mathematical Sciences, Shanghai Key Laboratory of PMMP, East China Normal University, Shanghai 200062,
People's Republic of China}
\email{zhiqiangwzy@163.com}

\author[K. Jiang]{Kan Jiang}
\address[K. Jiang]{Department of Mathematics, Ningbo University,
People's Republic of China}
\email{jiangkan@nbu.edu.cn}

\author[W. Li]{Wenxia Li}
\address[W. Li]{School of Mathematical Sciences, Shanghai Key Laboratory of PMMP, East China Normal University, Shanghai 200062,
People's Republic of China}
\email{wxli@math.ecnu.edu.cn}

\author[B. Zhao]{Bing Zhao}
\address[B. Zhao]{Department of Mathematics, Ningbo University,
People's Republic of China}
\email{zhaobing@nbu.edu.cn}

\date{\today}
\subjclass[2010]{Primary: 28A80, Secondary:11K55}

\begin{abstract}
 Let  $C$ be  the middle-third Cantor set. In this paper, we show that for every $x\in [0,4]$, there exist $x_1, x_2, x_3, x_4 \in C$ such that $x= x_1^2+x_2^2+x_3^2+x_4^2$, which was
  conjectured  in [Athreya, J. S.; Reznick, B.; Tyson, J. T. Cantor set arithmetic. \emph{Amer. Math. Monthly} 126 (2019), no. 1, 4–17].
\end{abstract}
\keywords{middle-third Cantor set, middle-$\frac{1}{\alpha}$ Cantor set, sum of four squares.}

\maketitle

\section{Introduction}

The middle-third Cantor set
\begin{equation*}\label{}
  C =\left\{ \sum_{i=1}^{\infty} \frac{\varepsilon_i}{3^i}: \varepsilon_i \in \{0,2\} \right\}
\end{equation*}
is a  classical object in fractal geometry.
The arithmetic on middle-third Cantor set has been studied in
 \cite{Athreya-Reznick-Tyson-2019,Cabrelli-Hare-Molter-1997,Jiang-Xi-2018,Gu-Jiang-Xi-Zhao2019,Mendes-Fernando-1994, Pawlowicz-2013,Utz-1951}.
The first classical result is that the set
\begin{equation}\label{substraction_Cantor_set}
  C-C:=\{x-y: x,y \in C\}
\end{equation}
equals to the interval $[-1,1]$.
The proof of (\ref{substraction_Cantor_set}) was first given by H. Steinhaus in 1917.
The result was rediscovered by J. F. Randolph in 1940 \cite{Randolph-1940}.
Using the symmetry of $C$, we can deduce that
\[ C+C=C+(1-C)=1+(C-C)=[0,2], \]
where $C+C:=\{x+y: x,y \in C\}$.
The multiplication and division on middle-third Cantor set were discussed in \cite{Athreya-Reznick-Tyson-2019}. They proved that
\[ \mathcal{L}(C\cdot C) \ge \frac{17}{21} \text{ and } \frac{C}{C}=\bigcup_{n=-\infty}^\infty \left[ \frac{2}{3}\cdot 3^n, \frac{3}{2} \cdot 3^n \right] \cup \{0\}, \]
where $C \cdot C:=\{xy: x,y \in C\}$, $\frac{C}{C}:=\left\{ \frac{x}{y}: x,y \in C, y \ne 0 \right\}$
and $\mathcal{L}$ denotes the Lebesgue measure on $\R$.
Gu, Jiang, Xi and Zhao \cite{Gu-Jiang-Xi-Zhao2019} gave the complete topological structure of $C \cdot C$.  Moreover, they also proved that the Lebesgue measure of $C \cdot C$ is about $0.80955$.

The main motivation of this paper is due to a conjecture posed by Athreya, Reznick and Tyson \cite{Athreya-Reznick-Tyson-2019}. They conjectured  $\{x_1^2+x_2^2+x_3^2+x_4^2:\; x_i\in C\}=[0,4]$ and
 claimed that there is strong numerical evidence supporting it. In this paper, we will prove this conjecture.

Fixing $\alpha >1$, let $C_\alpha $ (the middle-$\frac{1}{\alpha }$ Cantor set) be generated by the iterated function system $\Phi=\left\{ f_1(x)=rx, f_2(x)=rx+1-r \right\}$
with $r=\frac{1}{2}\left (1-\frac{1}{\alpha }\right )$. Thus the classical middle-third Cantor set $C=C_3$. In the present paper we prove

\begin{theorem}\label{general-theorem}
  Let $C_\alpha$ be the middle-$\frac{1}{\alpha}$ Cantor set for $\alpha>1$.
  Then
  \begin{equation*}\label{general-throrem-statement}
    \{x_1^2+x_2^2+x_3^2+x_4^2:\; x_i\in C_\alpha \}=[0,4] \textrm{ if and only if } \alpha \ge 3.
  \end{equation*}
\end{theorem}

This paper is organized as follows. In section 2, we  discuss the set $\{x_1^2+x_2^2+x_3^2:\; x_i\in C_\alpha \}$. The  proof of Theorem \ref{general-theorem} is arranged in the  section 3.

\section{Sum of three squares}

As stated in the previous section, $C_\alpha $ is the unique nonempty compact set satisfying
$$
C_\alpha =f_1(C_\alpha )\cup f_2(C_\alpha )=rC_\alpha \cup (rC_\alpha +1-r)
$$
where $r=\frac{1}{2}\left (1-\frac{1}{\alpha}\right )$. It follows that if $x\in C_\alpha$, then $rx \in C_\alpha$. We will use this simple observation in Lemma \ref{lemma-contract}.
For each positive integer $n$ let
$$
\mathcal{F}_n=\{f_\sigma ([0,1]):\; \sigma \in \{1,2\}^n\}\;\;\textrm{and}\;\; F_n=\bigcup _{A\in \mathcal{F}_n}A,
$$
where $f_\sigma (x)=f_{\sigma _1}\circ f_{\sigma _2}\circ \cdots \circ f_{\sigma _n}(x)$ for $\sigma =\sigma _1\sigma _2\cdots \sigma _n\in \{1,2\}^n$. Then the sequence $F_n,n=1,2,\cdots,$ of nonempty compact sets is decreasing and
$$
C_\alpha=\bigcap _{n=1}^\infty F_n=\bigcap _{n=1}^\infty \bigcup _{\sigma \in \{1,2\}^n}f_\sigma ([0,1]).
$$
It is easy to see that for $\sigma =\sigma _1\sigma _2\cdots \sigma _n\in \{1,2\}^n$
$$
f_\sigma (0)=\frac{1-r}{r}\sum _{k=1}^n(\sigma _k-1)r^k
$$
and  so
$$
f_\sigma ([0,1])=[f_\sigma (0), f_\sigma (1)]=\left [\frac{1-r}{r}\sum _{k=1}^n(\sigma _k-1)r^k,  \frac{1-r}{r}\sum _{k=1}^n(\sigma _k-1)r^k+r^n\right ].
$$

Each element of $\mathcal{F}_n$, called an $n$-level basic interval, has length $r^n$. For an $n$-level basic interval $f_\sigma ([0,1])$, it contains two $(n+1)$-level basic intervals $f_{\sigma 1}([0,1])$ and $f_{\sigma 2}([0,1])$. $f_\sigma ([0,1])$ shares the same left endpoint with $f_{\sigma 1}([0,1])$, and shares the same right endpoint with $f_{\sigma 2}([0,1])$. The length of the open interval $f_{\sigma }([0,1])\setminus (f_{\sigma 1}([0,1])\cup f_{\sigma 2}([0,1]))$
is $\frac{1}{\alpha }$  times that of   $f_{\sigma }([0,1])$.

Denote by $L_n$ the collection of left endpoints of all $n$-level basic intervals.
For $u \in L_n$, we associate $u$ with an $n$-level basic interval
\[ I_u= [ u,\; u+r^n ] \]
and two $(n+1)$-level basic intervals denoted by
\[ I_{u,0}=[u, u+r^{n+1}],\; I_{u,1}=[u+(1-r)r^n, u+r^n]. \]

The key to discuss the sum of squares of Cantor set is the following lemma, which is an easy exercise in real analysis.

\begin{lemma}\label{decreasing_compact_continuity}
  Let $\varphi: \R^d \to \R$ be continuous. If $\{ K_j\}_{j \in \N}$ is a decreasing sequence of nonempty compact subsets of $\R^d$, then \[ \varphi \left( \bigcap_{j=1}^\infty K_j \right)=\bigcap_{j=1}^\infty \varphi(K_j). \]
\end{lemma}
\begin{proof}
  Since $\bigcap_{j=1}^\infty K_j \subseteq K_n$ for every $n \in \N$, we have $\varphi \left( \bigcap_{j=1}^\infty K_j \right) \subseteq \bigcap_{j=1}^\infty \varphi(K_j)$.
  Conversely, assume that $y \in \bigcap_{j=1}^\infty \varphi(K_j)$.
  For every $j$, we can find $x_j \in K_j$ such that $\varphi(x_j)=y$.
  Since $K_1$ is compact, by  Bolzano–Weierstrass Theorem, there is a convergent subsequence $x_{n_j} \to x$.
  Since $\varphi$ is continuous, we have $\varphi(x)=y$.
  Note that the sequence $\{ x_{n_j} \}_{j \ge m}$ is in $K_m$ for every $m \in \N$.
  It follows from compactness that $x \in K_m$ for every $m \in \N$.
  Therefore, $y=\varphi(x) \in \varphi \left( \bigcap_{j=1}^\infty K_j \right)$, which completes the proof.
\end{proof}

Define  functions $g: \R^3 \to \R$ and  $f: \R^4 \to \R$ by letting
$$
g(x_1,x_2,x_3)= x_1^2+x_2^2+x_3^2
$$
and
$$
f(x_1,x_2,x_3, x_4)=g(x_1,x_2,x_3)+x_4^2= x_1^2+x_2^2+x_3^2+x_4^2.
$$
For a positive integer $k$ and a nonempty set $A\subseteq {\mathbb R}$, denote $A^k=\{(x_1, \cdots , x_k): x_i\in A\}$.
In order to show $f(C_\alpha^4)=[0,4]$, we need to discuss the set $g(C_\alpha^3)$ and find some intervals in $g(C_\alpha^3)$.
Note that $C_\alpha^3= \bigcap_{n=1}^\infty F_n^3$. Applying Lemma \ref{decreasing_compact_continuity} for the continuous function $g$, we obtain the following corollary.
\begin{corollary}\label{colollary-1}
  $g(C_\alpha^3)=\bigcap_{n=1}^\infty g(F_n^3)$.
\end{corollary}

If an interval $I \subseteq g(F_n^3)$ for every $n \in \N$, then $I \subseteq g(C_\alpha^3)$. The following two lemmas give a sufficient condition to find intervals in $g(C_\alpha^3)$.

\begin{lemma}\label{three-squares-overlap-lemma}
  Let $\alpha \ge 3$.
  For any $u,v,w \in L_n$, if
  \begin{equation}\label{con1}
    \max\{ u,v,w \} > 0
  \end{equation}
  and
  \begin{equation}\label{con2}
    4 (1-r) \max\{ u, v, w \} \le 2(u+v+w) + (1+2r)r^n ,
  \end{equation}
  then
  \[g(I_u\times I_v \times I_w )=g((I_{u,0}\cup I_{u,1})\times (I_{v,0}\cup I_{v,1})\times (I_{w,0}\cup I_{w,1})).\]
\end{lemma}
\begin{proof}
At first we have $r=\frac{1}{2}\left (1-\frac{1}{\alpha }\right )\in [1/3,  1/2)$ since $\alpha \geq 3$.
  Write $t=u^2+v^2+w^2$.
  Without loss of generality, we can assume that $ u \ge v \ge w $.
  By (\ref{con1}) we have $u > 0$ and so $u\geq f_{1^{n-1}2}(0)=(1-r)r^{n-1} > r^n$. In addition, (\ref{con2}) reduces to
  \begin{equation}\label{three-squares-overlap-condition}
    2v+2w+(1+2r)r^n \ge 2(1-2r)u.
  \end{equation}

  It is routine to verify that
  {\small
  \[ g(I_{u,1}\times I_{v,0} \times I_{w,0} )= \left[ t+2u(1-r)r^n+ (1-r)^2 r^{2n}, \; t+2(u+rv+rw)r^n+(1+2r^2)r^{2n} \right] , \]
  \[ g(I_{u,1}\times I_{v,0} \times I_{w,1} )=\left[ t+2(u+w)(1-r)r^n+ 2(1-r)^2 r^{2n},\; t+2(u+rv+w)r^n+(2+r^2)r^{2n} \right] , \]
  \[ g(I_{u,1}\times I_{v,1} \times I_{w,0} )= \left[ t+2(u+v)(1-r)r^n+ 2(1-r)^2 r^{2n},\; t+2(u+v+rw)r^n+(2+r^2)r^{2n} \right] , \]
  and
  \[ g(I_{u,1}\times I_{v,1} \times I_{w,1} )= \left[ t+2(u+v+w)(1-r)r^n+ 3(1-r)^2 r^{2n},\; t+2(u+v+w)r^n+3 r^{2n} \right] .\] }

Note that
\begin{equation*}
  \begin{split}
     & t+2(u+rv+rw)r^n+(1+2r^2)r^{2n} - (t+ 2(u+w)(1-r)r^n+2(1-r)^2r^{2n}) \\
     = & 2(ru+rv+2rw-w)r^n + (4r-1)r^{2n} \\
     \ge & 2(4r-1)wr^n + (4r-1)r^{2n} >  0,
  \end{split}
\end{equation*}
and
\begin{equation*}
  \begin{split}
     & t+2(u+rv+w)r^n+(2+r^2)r^{2n} - (t+2(u+v)(1-r)r^n+2(1-r)^2r^{2n}) \\
     = & 2(ru+2rv-v+w)r^n + (4-r)r^{2n+1} \\
     \ge & 2(3r-1)vr^n + (4-r)r^{2n+1} >  0,
  \end{split}
\end{equation*}
and
\begin{equation*}
  \begin{split}
     & t+2(u+v+rw)r^n+(2+r^2)r^{2n} - (t+ 2(u+v+w)(1-r)r^n+3(1-r)^2r^{2n}) \\
     = & 2(ru+rv+2rw-w)r^n + (6r-2r^2-1)r^{2n} \\
     \ge & 2(4r-1)wr^n + (6r-2r^2-1)r^{2n}
     >  0.
  \end{split}
\end{equation*}

  Therefore, we have
  \begin{equation}\label{right-part}
  \begin{split}
      & g(I_{u,1}\times ( I_{v,0} \cup I_{v,1}) \times ( I_{w,0}\cup I_{w,1} ) )\\
    = & \left[ t+2u(1-r)r^n+ (1-r)^2 r^{2n},\; t+2(u+v+w)r^n+3r^{2n} \right].
    \end{split}
  \end{equation}

  It is also routine to verify that
  {\small
  \[ g(I_{u,0}\times I_{v,0} \times I_{w,0})=\left[ t,\; t+2(u+v+w)r^{n+1}+3r^{2n+2} \right],\]
  \[ g(I_{u,0}\times I_{v,0} \times I_{w,1} )= \left[ t+2w(1-r)r^n+(1-r)^2 r^{2n},\; t+2(ru+rv+w)r^{n}+(1+2r^2) r^{2n} \right], \]
  \[ g(I_{u,0}\times I_{v,1} \times I_{w,0} )= \left[ t+2v(1-r)r^n+(1-r)^2 r^{2n},\; t+2(ru+v+rw)r^{n}+(1+2r^2) r^{2n} \right], \]
  and
  \[ g(I_{u,0}\times I_{v,1} \times I_{w,1} )= \left[ t+2(v+w)(1-r)r^n+2(1-r)^2 r^{2n},\; t+2(ru+v+w)r^{n}+(2+r^2) r^{2n} \right]. \] }
  Since $u > r^n $, we have
\begin{equation*}
  \begin{split}
     & t+2(u+v+w)r^{n+1}+3r^{2n+2} - (t+ 2w(1-r)r^n+(1-r)^2r^{2n}) \\
     = & 2(ru+rv+2rw-w)r^n + (2r^2+2r-1)r^{2n} \\
     \ge & 2(3r-1)wr^n + 2ur^{n+1} + (2r-1)r^{2n}\\
     > & 2(3r-1)wr^n  + (4r-1)r^{2n}
     > 0 ,
  \end{split}
\end{equation*}
and
\begin{equation*}
  \begin{split}
     & t+2(ru+rv+w)r^n+(1+2r^2)r^{2n} - (t+ 2v(1-r)r^n+(1-r)^2r^{2n}) \\
     = & 2(ru+2rv-v+w)r^n + (r+2)r^{2n+1} \\
     \ge & 2(3r-1)vr^n + (r+2)r^{2n+1}
     >  0 ,
  \end{split}
\end{equation*}
and
\begin{equation*}
  \begin{split}
     & t+2(ru+v+rw)r^n+(1+2r^2)r^{2n} - (t+ 2(v+w)(1-r)r^n+2(1-r)^2r^{2n}) \\
     = & 2(ru+rv+2rw-w)r^n + (4r-1)r^{2n} \\
     \ge & 2(4r-1)wr^n + (4r-1)r^{2n}
     > 0 .
  \end{split}
\end{equation*}
   Therefore, we have
  \begin{equation}\label{left-part}
    g(I_{u,0}\times ( I_{v,0} \cup I_{v,1}) \times ( I_{w,0}\cup I_{w,1} ) )= \left[ t, \; t+2(ru+v+w)r^{n}+(2+r^2) r^{2n} \right].
  \end{equation}

  It follows from condition (\ref{three-squares-overlap-condition}) that
  \[ t+2(ru+v+w)r^{n}+(2+r^2) r^{2n} \ge t+2u(1-r)r^n+ (1-r)^2 r^{2n}.\]
  Thus, the intervals in (\ref{right-part}) and (\ref{left-part}) overlap and so
  \[g( ( I_{u,0} \cup I_{u,1} ) \times ( I_{v,0} \cup I_{v,1}) \times ( I_{w,0}\cup I_{w,1} ) ) = \left[ t, \; t+2(u+v+w)r^n+3r^{2n} \right]. \]
  Note that
  \[ g(I_u\times I_v \times I_w )= \left[  t, \; t+2(u+v+w)r^n+3r^{2n} \right]. \]
  Therefore, we conclude that
  \[g(I_u\times I_v \times I_w )=g((I_{u,0}\cup I_{u,1})\times (I_{v,0}\cup I_{v,1})\times (I_{w,0}\cup I_{w,1})),\]
  as desired.
\end{proof}

\begin{lemma}\label{main-lemma}
  Let $\alpha \ge 3$. For any $u,v,w \in L_n$, if
  \begin{equation}\label{condition-interval}
    2(1-r) \max\{ u, v, w \} + (1-2r)r^n \le u+v+w ,
  \end{equation}
  then
  \[ g(I_u\times I_v \times I_w )\subseteq g(C_\alpha^3). \]
\end{lemma}
\begin{proof}
  For $k \ge n$, we define
  \[
    \mathcal{F}_{1,k}= \{ I \in \mathcal{F}_k: I \subseteq I_u \}, \;
    \mathcal{F}_{2,k}= \{ I \in \mathcal{F}_k: I \subseteq I_v \}, \;
    \mathcal{F}_{3,k}= \{ I \in \mathcal{F}_k: I \subseteq I_w \},
  \]
  and \[ F_{1,k}= \bigcup _{A\in \mathcal{F}_{1,k}}A,\;\; F_{2,k}= \bigcup _{A\in \mathcal{F}_{2,k}}A,\;\; F_{3,k}= \bigcup _{A\in \mathcal{F}_{3,k}}A. \]
  By Corollary \ref{colollary-1}, it suffices to show that for $k \ge n$,
  \begin{equation}\label{iqe}
   g(I_u\times I_v \times I_w)\subseteq g(F_{1,k} \times F_{1,k} \times F_{1,k}).
   \end{equation}
  We now prove it by induction on $k$.

  When $k=n$, we have $F_{1,n}=I_u, F_{2,n}=I_v, F_{3,n}=I_w$, and thus
  \[ g(I_u\times I_v \times I_w) \subseteq g(F_{1,n} \times F_{2,n} \times F_{3,n}). \]
  Next, assume that (\ref{iqe}) is true for some $m \ge n$, i.e.,
  \begin{equation}\label{gna}
   g(I_u\times I_v \times I_w) \subseteq g(F_{1,m} \times F_{2,m} \times F_{3,m}).
   \end{equation}
  Then, taking $x \in g(I_u\times I_v \times I_w)$,
  it follows from (\ref{gna}) that there exist $u',v',w' \in L_m$ such that
  \[ I_{u'} \subseteq I_{u},\; I_{v'} \subseteq I_{v},\; I_{w'} \subseteq I_{w} \;\;
  \textrm{and}\;\; x\in g(I_{u'} \times I_{v'} \times I_{w'}).\]

  Now condition (\ref{condition-interval}) implies that $\max\{ u,v,w \} > 0$, and thus
  \begin{equation*}\label{}
    \max\{ u',v',w' \} > 0.
  \end{equation*}
  Moreover, it follows from (\ref{condition-interval}) that
  \[
    \begin{split}
       (1-2r) \max\{ u',v',w' \}
       & \le (1-2r)\max\{ u,v,w \} +(1-2r)r^n \\
       & \le u+v+w - \max\{ u,v,w \} \\
       & \le u'+v'+w' - \max\{ u',v',w' \},
    \end{split}
  \]
  i.e.,
  \begin{equation*}\label{}
    2(1-r) \max\{ u', v', w' \} \le u'+v'+w'.
  \end{equation*}
  Thus, applying Lemma \ref{three-squares-overlap-lemma}, there exist $i,j,\ell \in \{0,1\}$ such that
  \[ x\in g(I_{u',i} \times I_{v',j} \times I_{w',\ell}).\]
  Obviously, we have $ I_{u',i} \in \mathcal{F}_{1,m+1},\; I_{v',j} \in \mathcal{F}_{1,m+1}$ and $ I_{w',\ell} \in \mathcal{F}_{1,m+1} $.
  Therefore,
  \[ x\in g(F_{1,m+1} \times F_{2,m+1} \times F_{3,m+1}). \]
  This shows that (\ref{iqe}) is true for $k=m+1$.
\end{proof}

\begin{corollary}\label{intervals-of-three-squares}
  For $\alpha \ge 3$,
  \[ \left[ a, b \right] \cup [2(1-r)^2,3] \subseteq g(C_\alpha^3), \]
  where $a=2r^4-4r^3+3r^2-2r+1$ and $b=r^4-2r^3+5r^2-2r+1$.
\end{corollary}
\begin{proof}
Note that
\begin{equation*}
\begin{split}
&g\left( \left[ 0, r \right] \times \left[ 1-r, 1 \right] \times \left[ 1-r, 1 \right]\right) \cup  g\left(\left[ 1-r, 1 \right] \times \left[ 1-r, 1 \right] \times \left[ 1-r, 1 \right]  \right) \\
&=\left[ 2(1-r)^2, 2+r^2 \right] \cup \left[ 3(1-r)^2, 3 \right]
   =[2(1-r)^2,3]
  \end{split}
  \end{equation*}
  and
  \[g\left( \left[ r-r^2, r \right] \times \left[ r-r^2, r \right] \times \left[ 1-r, 1-r+r^2 \right]  \right)= [a,b]. \]
We claim  that the intervals $ g\left( \left[ 0, r \right] \times \left[ 1-r, 1 \right] \times \left[ 1-r, 1 \right]  \right)$, $g\left( \left[ 1-r, 1 \right] \times \left[ 1-r, 1 \right] \times \left[ 1-r, 1 \right]  \right)$  and $g\left( \left[ r-r^2, r \right] \times \left[ r-r^2, r \right] \times \left[ 1-r, 1-r+r^2 \right]  \right)$ are all included in $g(C_\alpha^3)$. By Lemma \ref{main-lemma} these just are done by checking condition (\ref{condition-interval}). In fact, we have
\begin{equation*}
  2(1-r) \cdot (1-r) + (1-2r)r -2(1-r) = -r < 0,
\end{equation*}
\begin{equation*}
  2(1-r) \cdot (1-r) + (1-2r)r -3(1-r) = -1 < 0,
\end{equation*}
and
\begin{equation*}
  \begin{split}
     & 2(1-r) \cdot (1-r) + (1-2r)r^2 - (2(r-r^2)+(1-r)) \\
     = & -2r^3+5r^2-5r+1 = -r(2r-1)(r-2)-(3r-1)
     <  0.
  \end{split}
\end{equation*}
\end{proof}

\section{The proof of Theorem \ref{general-theorem}}

For $E \subseteq \R$ and $t \in \R$, we define $ t\cdot E=\{tx: x\in E \} $.

\begin{lemma}\label{lemma-contract}
  If $E \subseteq f(C_\alpha^4)$, then $r^2 \cdot E \subseteq f(C_\alpha^4)$.
  Similarly, if $E \subseteq g(C_\alpha^3)$, then $r^2 \cdot E \subseteq g(C_\alpha^3)$.
\end{lemma}
\begin{proof}
  Assume that $E \subseteq f(C_\alpha^4)$.
  For $x\in E$, there are $x_1, x_2, x_3, x_4 \in C_\alpha$ such that $x= x_1^2+x_2^2+x_3^2+x_4^2$.
  Then $r^2 x= (rx_1)^2+(rx_2)^2+(rx_3)^2+(rx_4)^2 \in f(C_\alpha^4)$.
  It follows that $r^2 \cdot E \subseteq f(C_\alpha^4)$.

  Similarly, the result for $g(C_\alpha^3)$ can be proved.
\end{proof}

\begin{lemma}\label{lemma-reducton}
  $f(C_\alpha^4)=[0,4]$ if and only if $(4r^2, 4] \subseteq f(C_\alpha^4)$.
\end{lemma}
\begin{proof}
  Note that
  \[ 0 \in f(C_\alpha^4) \textrm{ and } (0,4]=\bigcup_{n=0}^\infty r^{2n} \cdot \left( 4r^2, 4 \right]. \]
  The sufficiency follows from Lemma \ref{lemma-contract}.
\end{proof}

Now we are ready to prove Theorem \ref{general-theorem}.
\begin{proof}[The proof of Theorem \ref{general-theorem}]
  For $1 < \alpha < 3$, we have $0< r < \frac{1}{3}$, which implies $4r^2 < (1-r)^2$. Thus,
  \begin{equation*}\label{}
    (4r^2, (1-r)^2) \cap f(C_\alpha^4) \subseteq (4r^2, (1-r)^2) \cap f(F_1^4) = \emptyset.
  \end{equation*}
  Therefore, it suffices to show   $f(C_\alpha^4)=[0,4]$ when $\alpha \ge 3$.

  Assume that $\alpha \ge 3$. Note that $\frac{1}{3} \le r < \frac{1}{2}$.
  Then we have $(1-r)^2 \le 4r^2$.
  By Lemma \ref{lemma-reducton}, it suffices to prove that
  \begin{equation}\label{goal}
  ( (1-r)^2, 4] \subseteq f(C_\alpha^4).
  \end{equation}

  In Corollary \ref{intervals-of-three-squares}, we have $[2(1-r)^2,3] \subseteq g(C_\alpha^3)$. Thus
  \begin{equation}\label{one}
   f(C_\alpha^4) \supseteq  f(C_\alpha ^3\times \{ 0,1 \})=(g(C_\alpha^3)+0^2)\cup (g(C_\alpha^3)+1^2)\supseteq
   \left[ 2(1-r)^2 , 4 \right]
  \end{equation}

  Applying Corollary \ref{intervals-of-three-squares} and Lemma \ref{lemma-contract}, we have
  \[  g(C_\alpha^3)\supseteq [ a r^{2n}, b r^{2n} ] \cup [ 2(1-r)^2 r^{2n}, 3\cdot r^{2n} ]\;\;\textrm{for}\;\; n=0,1,2,\cdots . \]
  where $a,b$ are given in Corollary \ref{intervals-of-three-squares}.
    Note that for each positive integer $n$
  \begin{equation}\label{partial-result-1}
  \begin{split}
  f(C_\alpha^4)&\supseteq f\left( C_\alpha ^3  \times \left\{ 1-r, 1-r+r^{2n} \right\}\right)\\
    &=\left (g(C_\alpha^3)+(1-r)^2\right )\cup \left( g(C_\alpha^3)+(1-r+r^{2n})^2\right )\\
    &\supseteq [a r^{2n-2}+(1-r)^2, b r^{2n-2}+(1-r)^2]\\
    &\;\;\;\; \cup [a r^{2n-2}+(1-r+r^{2n})^2, b r^{2n-2}+(1-r+r^{2n})^2]\\
    &=[a r^{2n-2}+(1-r)^2, b r^{2n-2}+(1-r+r^{2n})^2]\\
    &\supseteq [ (1-r)^2 + a r^{2n-2} , (1-r)^2+(b+2r^2-2r^3)r^{2n-2} ]
  \end{split}
  \end{equation}
  where the last equality and the last inclusion hold because
  \begin{equation*}
  \begin{split}
     & br^{2n-2} + (1-r)^2 - (ar^{2n-2}+(1-r+r^{2n})^2 )\\
     = & (b-a)r^{2n-2} -2(1-r)r^{2n}-r^{4n} \\
     = & (4r-r^2-r^{2n})r^{2n}
     \ge (4r-2r^2)r^{2n}>  0,
  \end{split}
\end{equation*}
and
\begin{equation*}
  br^{2n-2}+(1-r+r^{2n})^2 - ((1-r)^2+(b+2r^2-2r^3)r^{2n-2} )= r^{4n} > 0.
\end{equation*}

 Note that for each positive integer $n$
  \begin{equation}\label{partial-result-2}
  \begin{split}
  f(C_\alpha^4)&\supseteq f\left( C_\alpha ^3  \times \left\{ 1-r, 1-r+r^{2n},1-r+r^{2n-1}-r^{2n} \right\}\right)\\
    &=\left (g(C_\alpha^3)+(1-r)^2\right )\cup \left( g(C_\alpha^3)+(1-r+r^{2n})^2\right )\\
    &\;\;\;\;\;\cup \left( g(C_\alpha^3)+(1-r+r^{2n-1}-r^{2n})^2\right )\\
    &\supseteq [ 2(1-r)^2 r^{2n}+(1-r)^2, 3r^{2n}+(1-r)^2]\\
   &\;\;\;\; \cup [ 2(1-r)^2 r^{2n}+(1-r+r^{2n})^2, 3r^{2n}+(1-r+r^{2n})^2]\\
   &\;\;\;\; \cup [ 2(1-r)^2 r^{2n}+(1-r+r^{2n-1}-r^{2n})^2, 3r^{2n}+(1-r+r^{2n-1}-r^{2n})^2]\\
   &=[ 2(1-r)^2 r^{2n}+(1-r)^2, 3r^{2n}+(1-r+r^{2n-1}-r^{2n})^2]\\
   & \supseteq  [ (1-r)^2 + 2(1-r)^2 r^{2n}, (1-r)^2+(2-r+2r^2)r^{2n-1} ]
  \end{split}
  \end{equation}
  where the last equality and inclusion hold because
  \begin{equation*}
  \begin{split}
     & 3r^{2n}+(1-r)^2 - (2(1-r)^2r^{2n}+(1-r+r^{2n})^2) \\
     = & (6r-2r^2-1-r^{2n})r^{2n} \\
     \ge & (6r-3r^2-1)r^{2n}= [3(1-r)r+ 3r-1]r^{2n}
     >  0,\\
    & 3r^{2n} + (1-r+r^{2n})^2 - (2(1-r)^2r^{2n} +( 1-r+r^{2n-1}-r^{2n})^2) \\
    = & -2r^{2n-1}+7r^{2n}-2r^{2n+2}-r^{4n-2}+2r^{4n-1} \\
    = & 2(3r-1)r^{2n-1}+(1-2r^2)r^{2n}-r^{4n-2}+2r^{4n-1} \\
    > & 2(3r-1)r^{2n-1}+r^{2n+1}-r^{4n-2}+2r^{4n-1}  \\
    \ge & 2(3r-1)r^{2n-1}-r^{4n-2}+3r^{4n-1} = 2(3r-1)r^{2n-1}+(3r-1)r^{4n-2}
    \ge  0,
  \end{split}
\end{equation*}
and
\begin{equation*}
    3r^{2n} + ( 1-r+r^{2n-1}-r^{2n})^2 - ((1-r)^2 + (2-r+2r^2)r^{2n-1}) = (r^{2n-1}-r^{2n})^2>0.
\end{equation*}

Note that
\begin{equation*}
\begin{split}
&a-(2-r+2r^2)r = 2r^4-6r^3+4r^2-4r+1 \\
&=(2r-1)r^3-(5r^2-4r+1)r -(3r-1) <0,
\end{split}
  \end{equation*}
  and
  \begin{equation*}
\begin{split}
  &(b+2r^2-2r^3)-2(1-r)^2=r^4-4r^3+5r^2+2r-1\\
  & = r^4+2r^2(1-2r)+(3r-1)(r+1) >0.
  \end{split}
  \end{equation*}

 Combining the above inequalities, (\ref{partial-result-1}) and (\ref{partial-result-2}) we have that for every $n \ge 1$,
  \begin{equation*}\label{}
   f(C_\alpha^4)\supseteq  [ (1-r)^2+2(1-r)^2 r^{2n}, (1-r)^2+2(1-r)^2 r^{2n-2} ].
  \end{equation*}
   Therefore,
   \begin{equation}\label{two}
  f(C_\alpha^4)\supseteq \bigcup_{n=1}^\infty[(1-r)^2(1+2r^{2n}), (1-r)^2(1+2r^{2n-2})] = ((1-r)^2, 3(1-r)^2].
  \end{equation}
  By (\ref{one}) and (\ref{two}), we have
  $ f(C_\alpha^4)\supseteq  ((1-r)^2, 3(1-r)^2] \cup [2(1-r)^2, 4]=( (1-r)^2, 4],$
  obtaining (\ref{goal}).
\end{proof}

\section*{Acknowledgements}
The second author was supported by NSFC No. 11701302 and K.C. Wong Magna Fund in Ningbo University. The second author was also supported by  Zhejiang Provincial Natural Science Foundation of China with No.LY20A010009.
The third author was supported by NSFC No. 11671147, 11971097 and Science and Technology Commission of Shanghai Municipality (STCSM) No. 13dz2260400.

\end{document}